\documentclass[11pt,a4wide]{article}
\usepackage{amsmath}
\usepackage{amssymb}
\usepackage{dsfont}
\usepackage{bbm}
\usepackage{graphicx}
\usepackage{psfrag}
\usepackage[amsmath,thmmarks]{ntheorem}
\usepackage{hyperref}
\usepackage{color}

\theoremstyle{plain}
\theoremheaderfont{\normalfont\bfseries}
\theorembodyfont{\slshape}
\theoremseparator{.}
\newtheorem{Theorem}{Theorem}[section]

\theoremstyle{plain}
\theoremindent0.0cm \theoremnumbering{greek}

\theoremstyle{plain} 
\theoremindent0.0cm \theoremnumbering{greek}
\newtheorem{proposition}[Theorem]{Proposition}

\theoremsymbol{}
\theoremnumbering{arabic}
\newtheorem{corollary}[Theorem]{Corollary}

\theoremstyle{plain} \theorembodyfont{\upshape} \theoremsymbol{\ensuremath{\Box}} \theoremseparator{}
\newtheorem{example}[Theorem]{Example}

\theoremstyle{plain}
\theoremsymbol{\ensuremath{\Box}}
\theoremseparator{.}
\newtheorem{definition}[Theorem]{Definition}

\theoremstyle{plain} \theoremsymbol{\ensuremath{\Box}} \theoremseparator{.} 
\newtheorem{remark}[Theorem]{Remark}

\theoremheaderfont{\sc} \theorembodyfont{\upshape} \theoremstyle{nonumberplain}
\theoremsymbol{\ensuremath{\Box}} 
\theoremseparator{:}

\theoremstyle{plain} \theoremsymbol{\ensuremath{\Box}} \theoremindent0.0cm \theoremnumbering{greek}

\title{$C_0$-semigroups for hyperbolic partial differential equations on a one-dimensional spatial domain}
\author{Birgit Jacob\footnote{ University of Wuppertal, Fachbereich C - Mathematik und Naturwissen\-schaften,
Arbeitsgruppe Funktionalanalysis, Gau\ss stra\ss e 20
D-42119 Wuppertal, Germany, bjacob@uni-wuppertal.de}
 \and Kirsten Morris\footnote{Department of Applied Mathematics,
University of Waterloo,
Waterloo, Ontario, Canada N2L 3G1, kmorris@uwaterloo.ca}
 \and Hans Zwart\footnote{University of Twente, Faculty of Electtical Engineering, Mathematics and Computer Science, Department of Applied Mathematics, P.O. Box 217, 7500 AE Enschede, The Netherlands, h.j.zwart@math.utwente.nl}}
\date{}
\begin{document}
\maketitle

\begin{abstract}
  Hyperbolic partial differential equations on a one-dimensional
  spatial domain are studied. This class of systems includes models of
  beams and waves as well as the transport equation and networks of
  non-homogeneous transmission lines. The main result of this paper is
  a simple test for $C_0$-semigroup generation in terms of the
  boundary conditions. The result is illustrated with several
  examples.
\end{abstract}


\section{Introduction and main result}

Consider the following class of partial differential
  equations 
\begin{align}\label{eqn:pde}
\frac{\partial x}{\partial t}(\zeta,t) =& \left( P_1
  \frac{\partial}{\partial \zeta} + P_0\right) ({\cal H}({\zeta})
x(\zeta,t)),\qquad \zeta \in [0,1], t\ge 0,\\
x(\zeta,0) =& x_0(\zeta),\nonumber
\end{align}
where $P_1$ is an invertible $n\times n$ Hermitian matrix, $P_0$ is a $n\times n$  matrix, 
${\cal H}(\zeta)$ is a positive $n\times n$ Hermitian matrix for a.e.~$\zeta\in (0,1)$  satisfying ${\cal H}, {\cal H}^{-1}\in L^{\infty}(0,1;\mathbb C^{n\times n})$. This class of Cauchy problems covers in particular the wave equation, the transport equation and the Timoshenko beam equation, and also coupled beam and wave equations.
These Cauchy problems are also known as Hamiltonian partial differential equations or port-Hamiltonian systems, see  \cite{JZ} ,\cite{VanDerSchaftMaschke_2002} and in particular the Ph.D thesis \cite{Villegas_2007}. 
The boundary conditions are of the form
\begin{equation}
\label{eq:bc}
   \tilde{W}_B\left[ \begin{smallmatrix} ({\cal H}x)(1,t) \\ ({\cal
         H}x)(0,t)\end{smallmatrix}\right]  =0, 
\end{equation}
where $\tilde{W}_B$ is an $n\times 2n$-matrix. 
Define
\begin{equation}\label{operatorA}
Ax:= \left( P_1 \frac{d}{d\zeta} + P_0\right) (x), \qquad x\in D(A), 
\end{equation}
 on $X_p:=L^p(0,1;\mathbb C^n)$, $1\le p <\infty$,  with the domain
\begin{equation}
\label{domainA}
D(A) := \left\{ x\in {\mathcal W}^{1,p}(0,1;\mathbb C^n)\mid  \tilde{W}_B\left[ \begin{smallmatrix} x(1) \\ x(0)\end{smallmatrix}\right] = 0 \right\}.
\end{equation}
Then the partial differential equation (\ref{eqn:pde}) with the
boundary conditions (\ref{eq:bc}) can be written as the abstract
differential equation
\[
   \dot{x} (t) = A {\cal H} x(t),\qquad   x(0)= x_0 .
\]

If we equip $X_2$ with the energy norm $\langle \cdot,{\cal H} \cdot\rangle$, then $A{\cal H}$ generates a contraction semigroup (or
an unitary $C_0$-group) on $(X_2,\langle \cdot,{\cal H} \cdot\rangle)$ if and only if $A$ is dissipative on $(X_2,\langle \cdot, \cdot\rangle)$(or $A$ and $-A$ are dissipative on $(X_2,\langle \cdot, \cdot\rangle)$, respectively)  \cite{AJ,JZ,GZM}. 
Matrix conditions to guarantee generation of a contraction semigroup
or of a unitary group have been obtained  \cite{AJ,JZ,GZM}. 
The following theorem extends these results.
\begin{Theorem}%
\label{theo1}  
 Let $W_B:=\tilde{W}_B\left[\begin{smallmatrix}  P_1 & - P_1 \\ I & I  \end{smallmatrix}\right]^{-1}$ and $\Sigma := \left[\begin{smallmatrix} 0 & I \\ I & 0 \end{smallmatrix}\right]$. 
\begin{enumerate}
\item The following statements are equivalent:
\begin{enumerate}
\item $A{\cal H}$ with domain $D(A{\cal H}):= \{x\in X_2\mid {\cal H}x\in D(A)\}={\cal H}^{-1}D(A)$ generates a contraction semigroup on $(X_2,\langle \cdot,{\cal H} \cdot\rangle)$;
\item Re$\,\langle Ax,x\rangle\le 0$ for every $x\in D(A)$;
\item Re$\, P_0\le 0$ and $u^*P_1 u- y^*P_1 y\le 0$  for every $\left[\begin{smallmatrix} u \\ y \end{smallmatrix}\right]\in \mbox{ ker}\, \tilde W_B$;
\item Re$\, P_0\le 0$, $W_B\Sigma W_B^*\ge 0$ and rank$\,\tilde{W}_B=n$.
\end{enumerate}
\item The following statements are equivalent:
\begin{enumerate}
\item $A{\cal H}$ with domain $D(A{\cal H}):= \{x\in X_2\mid {\cal H}x\in D(A)\}={\cal H}^{-1}D(A)$ generates a unitary $C_0$-group on $(X_2,\langle \cdot,{\cal H} \cdot\rangle)$;
\item Re$\,\langle Ax,x\rangle= 0$ for every $x\in D(A)$;
\item Re$\, P_0= 0$ and $u^*P_1 u- y^*P_1 y= 0$  for every $\left[\begin{smallmatrix} u \\ y \end{smallmatrix}\right]\in \mbox{ ker}\, \tilde W_B$;
\item Re$\, P_0= 0$, $W_B\Sigma W_B^*= 0$ and rank$\,\tilde{W}_B=n$.
\end{enumerate}
\end{enumerate}
\end{Theorem}

Theorem \ref{theo1} was proved  in
\cite[Theorem 7.2.4]{JZ} with the additional assumptions that $P_0^*=-P_0$ and rank$\,\tilde{W}_B=n$. The extension to non skew-adjoint matrices $P_0$ is in \cite{AJ}. However,  the equivalence with (c)  is not explicitly
shown in the above references and it is assumed that  rank$\,\tilde{W}_B=n$.  A short proof of Theorem \ref{theo1} is in
the following section. 

By the assumptions on ${\mathcal H}$ it is clear that the norm on
$(X_2,\langle \cdot,{\cal H} \cdot\rangle)$ is equivalent to the
standard norm on $X_2$. Hence if $A{\mathcal H}$ generates a contraction (or a
unitary group) with respect to the energy norm for some ${\mathcal H}$,
then it will generate a $C_0$-semigroup ($C_0$-group) on $X_2$ equipped with the standard norm as
well. 

The following corollary follows immediately.
\begin{corollary} \label{cor}  The following statements are equivalent:
\begin{enumerate}
\item $A$ generates a contraction semigroup on $(X_2,\langle\cdot , \cdot\rangle)$,
\item $A{\cal H}$ generates a contraction semigroup on $(X_2,\langle \cdot,{\cal H} \cdot\rangle)$.
\end{enumerate}
\end{corollary}

Corollary \ref{cor} implies that whether or not $A{\cal H}$ generates a contraction semigroup on the energy space $(X_2,\langle \cdot,{\cal H} \cdot\rangle)$ is independent of the Hamiltonian density ${\cal H}$: $A$ is the generator of a  contraction semigroup on $(X_2,\langle
\cdot, \cdot\rangle)$ if and only if $A{\cal H}$ generates a contraction semigroup on $(X_2,\langle
\cdot,{\cal H} \cdot\rangle)$.  The condition of a contraction
semigroup is essential here. For a counterexample, see Example \ref{ex2} or \cite[Section 6]{ZGMV}.
%
%
%
\begin{definition}
An operator ${\cal A}$ generates a quasi-contractive semigroup  if ${\cal A}- \omega I $ generates a contraction semigroup for some $\omega \in \mathbb R$.
\end{definition}

\begin{corollary}\label{corquasi}
If Re$\,P_0\le 0$ then  $A{\cal H}$ generates a quasi-contractive  semi\-group  on $(X_2,\langle\cdot ,{\cal H} \cdot\rangle)$ if and only if   $A{\cal H}$ generates a  contraction semigroup  on $(X_2,\langle\cdot ,{\cal H} \cdot\rangle)$.
\end{corollary}

The proof of Corollary \ref{corquasi} will be given in Section 2.

Theorem \ref{theo1} characterizes boundary conditions for which $A{\cal H}$ generates a 
contraction semigroup or a unitary group. However, other
boundary conditions may still lead to a $C_0$-semigroup. To
characterize those we diagonalize $P_1{\mathcal H}(\zeta)$. 
It is easy to see that the eigenvalues of $P_1{\mathcal H}(\zeta)$
are the same as the eigenvalues of ${\mathcal H}(\zeta)^{\frac{1}{2}}
  P_1{\mathcal H}(\zeta)^{\frac{1}{2}}$. Hence by Sylvester's Law of Inertia the number of positive and negative eigenvalues of
    $P_1{\mathcal H}(\zeta)$ equal those of $P_1$. We denote by $n_1$
    the number of positive and by $n_2=n-n_1$ the number of negative eigenvalues of
    $P_1$. Hence we can find
    matrices  such that
\begin{equation}\label{eqndiagonal}
 P_1{\cal H}(\zeta) =S^{-1}(\zeta)\begin{bmatrix} \Lambda(\zeta) & 0 \\ 0 & \Theta(\zeta)
\end{bmatrix}S(\zeta),\qquad \mbox{ a.e.~} \zeta\in(0,1),
\end{equation}
with $\Lambda(\zeta)$ and $\Theta(\zeta)$ diagonal matrices of size
$n_1 \times n_1$ and $n_2 \times n_2$, respectively. 

The main result  of this paper is the following theorem that  provides  easily checked  conditions for when the operator $A{\cal H}$ generates a $C_0$-semigroup on $X_p$. These cover the situation where $A{\cal H}$ may not generate a contraction semigroup.
\begin{Theorem}\label{Thmmain}
Assume that $S$, $\Lambda$ and $\Theta$ in
  (\ref{eqndiagonal}) are continuously differentiable  on $[0,1]$ and that rank$\, \tilde{W}_B=n$.
Define $Z^+(\zeta)$ to be the span of eigenvectors of $P_1 {\mathcal
  H}(\zeta)$ corresponding to its positive eigenvalues. Similarly, we
define $Z^-(\zeta)$  to be the span of eigenvectors of $P_1 {\mathcal
  H}(\zeta)$ corresponding to its negative eigenvalues. We write
$\tilde{W}_B$ as 
\begin{equation}
  \tilde{W}_B= \left[ \begin{matrix} W_{1} & W_{0} \end{matrix} \right]
\end{equation}
with $W_{1}, W_{0} \in {\mathbb C}^{n \times n}$. Then the following statements are equivalent:
\begin{enumerate}
\item The
 operator $A{\cal H}$ defined by (\ref{operatorA})--(\ref{domainA}) generates a $C_0$-semigroup on $X_p$. 
\item $ W_1 {\mathcal H}(1)  Z^+(1) \oplus W_0 {\mathcal H}(0) Z^-(0) = {\mathbb C}^n$.
  \end{enumerate}
\end{Theorem}

The proof of Theorem \ref{Thmmain} will be given in the next
section. 

\begin{remark}
\begin{enumerate}
\item In Kato \cite[ Chapter II]{K}, conditions on $P_1{\cal H}$ are given guaranteeing that $S$, $\Lambda$ and $\Theta$ are continuously differentiable. 
\item  In \cite{E},  a more restrictive version of  Theorem \ref{Thmmain} that applies when ${\cal H}=I$ and $p=2$ was proven by  a different approach. In \cite{E} estimates for the growth bound are given.  
\item Theorem \ref{Thmmain} implies that if $A{\cal H}$ generates a $C_0$-semigroup on one $X_p$, then $A{\cal H}$ generates  a $C_0$-semigroup on every $X_p$, $1\le p <\infty$. A similar statement 
does not hold for contraction semigroups. Example \ref{ex3}, given later in this paper, illustrates this point. 
\end{enumerate}
\end{remark}

%

\section{Proof of Theorems \ref{theo1} and \ref{Thmmain} and Corollary \ref{corquasi}}

\medskip

\noindent
{\bf Proof of Theorem \ref{theo1}:}

 Since
  the proof of Part 2 is similar to that of Part 1 we only present
  the details for Part 1.
  
The implication (a) $\Rightarrow$ (b) follows directly from the
Lumer-Phillips theorem and Lemma 7.2.3 in \cite{JZ}. 
Next we show the implication (b) $\Rightarrow$ (c).  It is easy to see that 
\begin{equation}
\label{eq:1HZ}
  \mathrm{Re}\langle Ax,x\rangle = x(1)^*P_1 x(1) - x(0)^*P_1 x(0) + \mathrm{Re} \int_0^1 x(\zeta)^* P_0 x(\zeta) d\zeta
\end{equation} 
holds for every $x\in D(A)$.
Choosing $x \in W^{1,2}(0,1;{\mathbb C}^n)$ with $x(0)=x(1)=0$, we obtain  Re$\, P_0\le 0$.  For every $u,y \in {\mathbb C}^n$ and every $\varepsilon >0$ there exists a function in $x \in W^{1,2}(0,1;{\mathbb C}^n)$ such that $x(0)=u$, $x(1)=y$ and the $L^2$-norm of $x$ is less than $\varepsilon$. Choosing this function in equation (\ref{eq:1HZ}) and letting $\varepsilon$ go to zero implies the second assertion in (c), see also Lemma 2.4 of \cite{AJ}.
The implication (d)
$\Rightarrow$ (a) follows from Theorem 2.3 of \cite{AJ}, see also \cite{GZM}.
Hence it remains to show (c) $\Rightarrow$ (d). 

We introduce the notation $f_1=x(1)$ and $f_0=
 x(0)$. Then the condition in (c) can be written as
\begin{equation}
  \label{eq:3}
  \left[ \begin{array}{cc} f_1^* & f_0^* \end{array} \right]
   \left[ \begin{array}{cc} P_1 & 0 \\ 0 & - P_1 \end{array} \right]
   \left[ \begin{array}{c} f_1\\ f_0 \end{array} \right] \leq 0, \quad \mbox{
     for } \left[ \begin{array}{c} f_1\\ f_0 \end{array} \right] \in
   \ker \tilde{W}_B.
\end{equation}
Since $\tilde{W}_B$ is an $n\times 2 n$ matrix, its kernel has
dimension $2n$ minus its rank. Hence this dimension will be larger or
equal to $n$. Since $P_1$ is an invertible Hermitian $n \times n$ matrix, the
matrix $\left[ \begin{smallmatrix} P_1 & 0 \\ 0 & -
    P_1 \end{smallmatrix} \right]$ will have $n$ positive and $n$
negative eigenvalues. This implies that if $v^*\left[ \begin{smallmatrix} P_1 & 0 \\ 0 & -
    P_1 \end{smallmatrix} \right]v \leq 0$ for all $v$ in a linear
subspace $V$, then $V$ has at most dimension $n$. Combining these two
facts,  the dimension of the kernel of $\tilde{W}_B$ equals
$n$, and so $\tilde{W}_B$ is a matrix of rank $n$. 

Defining
$\left[ \begin{smallmatrix} y_1\\ y_0 \end{smallmatrix} \right] =
\left[ \begin{smallmatrix} P_1 & -P_1 \\ I & I \end{smallmatrix}
\right] \left[ \begin{smallmatrix} f_1\\ f_0 \end{smallmatrix}
\right]$, 
and using (\ref{eq:3}), an easy calculation shows
\begin{equation}
  \label{eq:4}
  y_1^*y_0 + y_0^* y_1 \leq 0, \quad \mbox{
     for } \left[ \begin{array}{c} y_1\\ y_0 \end{array} \right] \in
   \ker W_B.
\end{equation}
We write $W_B$ as $W_B= [ W_1 \; W_2]$. Now  it is easy to see that  $W_1+W_2$ is invertible (we refer to page 87 in \cite{JZ} for the details). Defining $V:= (W_1+W_2)^{-1}(W_1-W_2)$, we obtain
\[ W_B = \frac{1}{2}(W_1+W_2)\left[ I+V, I-V\right].\] 
Let $\left[\begin{smallmatrix}  f \\
    e \end{smallmatrix} \right] \in \ker W_B$ be
arbitrary. By \cite[Lemma 7.3.2]{JZ}  there exists a vector $\ell$ such that $\left[\begin{smallmatrix}  f\\ e\end{smallmatrix} \right]=\left[\begin{smallmatrix}  I-V \\ -I-V \end{smallmatrix} \right]\ell$. This implies
\begin{equation}
\label{eq:1.11} 
  0 \ge  f^* e + e^* f= \ell^*(-2I+2V^*V)\ell,
\end{equation}
This inequality holds for any $\left[\begin{smallmatrix}  f \\ e \end{smallmatrix} \right] \in \ker W_B$. Since the $n\times 2n$ matrix $W_B$ has rank $n$, its kernel has dimension $n$, and so the set of vectors $\ell$ satisfying $\left[\begin{smallmatrix}  f \\ e \end{smallmatrix} \right]=\left[\begin{smallmatrix}  I-V \\ -I-V \end{smallmatrix} \right]\ell$ for some $\left[\begin{smallmatrix}  f \\ e \end{smallmatrix} \right] \in \ker W_B$ equals the whole space ${\mathbb K}^n$.  Thus (\ref{eq:1.11}) implies that $V^*V\le I$, and by \cite[Lemma 7.3.1]{JZ} we obtain $W_B\Sigma W_B^*\ge 0$.
\hfill${\Box}$
\bigskip

\noindent
{\bf Proof of Corollary \ref{corquasi}:}
 As $A{\mathcal H}-\omega I$ generates a contraction semigroup, Theorem \ref{theo1} implies $W_B\Sigma W_B^*\le  0$ and rank$\,\tilde{W}_B=n$. Thanks to Re$\, P_0\le 0$ and Theorem \ref{theo1}, finally $A{\mathcal H}$ generates a contraction semigroup.
\hfill${\Box}$
\bigskip

The following proposition is needed for the proof of Theorem \ref{Thmmain}.
\begin{proposition} \label{lemdiag}
(\cite[Theorem 3.3]{ZGMV} \cite[Theorem 13.3.1]{JZ} for $p=2$ and \cite[Theorem 3.3 and Section 7]{ZGMV} for $1\le p <\infty$)
Suppose $K, Q\in \mathbb C^{n\times n}$,  $\Lambda \in C^1([0,1];\mathbb C^{n_1\times n_1})$ is a diagonal real matrix-valued function with (strictly) positive functions on the diagonal and $\Theta\in C^1([0,1];\mathbb C^{n_2\times n_2})$, $n_1+n_2=n$, is a diagonal real matrix-valued function with (strictly) negative functions on the diagonal. We split a function
$g\in L^p(0,1;\mathbb C^n)$ as 
\begin{equation}
g(\zeta) =\begin{bmatrix} g_+(\zeta) \\ g_-(\zeta)\end{bmatrix},
\label{splitg}
\end{equation}
where $g_+(\zeta)\in \mathbb C^{n_1}$ and  $g_-(\zeta)\in \mathbb C^{n_2}$.

Then the operator $\tilde A: D(\tilde A)\subset X_p\rightarrow X_p$ defined by
\begin{align}
      \label{eq:1}
      \tilde A \begin{bmatrix} g_+\\ g_-\end{bmatrix}
      &=
      \frac{d}{d \zeta} \left(\begin{bmatrix} \Lambda & 0 \\ 0 & \Theta \end{bmatrix} \begin{bmatrix} g_+ \\ g_- \end{bmatrix}\right)
  \end{align}
\begin{align}
        D(\tilde A) &= \left\{ \begin{bmatrix} g_+ \\
        g_-\end{bmatrix} \! \in\!  W^{1,p}(0,1, {\mathbb C^n}) \mid
      \label{eq:2}
        K  \!\begin{bmatrix} \Lambda(1) g_+(1) \\ \Theta(0)
        g_-(0)\end{bmatrix}
        \!+ \!Q  \begin{bmatrix} \Lambda(0) g_+(0) \\ \Theta(1)
        g_-(1)\end{bmatrix} =0 \right\}
    \end{align}
 generates a $C_0$-semigroup on $X_p$ if and only if $K$ is invertible.
\end{proposition}

\noindent
{\bf Proof of Theorem \ref{Thmmain}:}
We define the new state variable $g := S x$. Since $S$ defines a boundedly invertible operator  on $L^p(0,1;{\mathbb C}^n)$, the operator $A{\mathcal H}$ generates a $C_0$-semigroup if and only if $S A{\mathcal H} S^{-1}$ generates a $C_0$-semigroup. We define
\[ \Delta := \begin{bmatrix} \Lambda & 0 \\ 0 & \Theta \end{bmatrix}.\]
Then the operator 
\begin{eqnarray}
  \nonumber
  (S A{\mathcal H} S^{-1} g)(\zeta) &=& \frac{d}{d \zeta}( \Delta (\zeta) g (\zeta)) +S(\zeta) \frac{d S^{-1}}{d\zeta}(\zeta)  \Delta(\zeta)  g (\zeta) \\
  \label{eq:8}
   && + S(\zeta) P_0{\mathcal H}(\zeta)S^{-1}(\zeta) g(\zeta)\\
\nonumber
  D(SA{\mathcal H}S^{-1}) &=& \{ g \in W^{1,p}(0,1; {\mathbb C^n}) \mid
  \tilde{W}_B\begin{bmatrix} ({\mathcal H}S^{-1} g)(1) \\  ({\mathcal H}S^{-1}g)(0) \end{bmatrix} =0\}.
\end{eqnarray}
Since the last two operators in (\ref{eq:8}) are bounded, $SA{\mathcal H}S^{-1}$ generates a $C_0$-semigroup if and only if the  operator  
\begin{align}
  \label{eq:10}
  A_S g &=  \frac{d}{d \zeta}( \Delta g ) \\
  \label{eq:11}
  D(A_S)  &= \left\{ g\in  W^{1,p}(0,1;\mathbb C^{n\times n}) \mid \tilde{W}_B\begin{bmatrix} ({\mathcal H}S^{-1} g)(1) \\  ({\mathcal H}S^{-1}g)(0) \end{bmatrix} =0\right\}
\end{align}
generates a $C_0$-semigroup on $X_p$. 
We split the matrices $W_{1} ({\mathcal H}S^{-1})(1)$ and $W_{0} ({\mathcal H}S^{-1})(0)$ as
\[  W_{1}({\mathcal H} S^{-1})(1)= \begin{bmatrix}  V_1 & V_2\end{bmatrix} \qquad  W_{0} ({\mathcal H}S^{-1})(0)= \begin{bmatrix}  U_1 & U_2 \end{bmatrix}, \]
where $U_1, V_1\in \mathbb C^{n\times n_1}$ and $U_2, V_2\in \mathbb C^{n\times n_2}$, and as in (\ref{split}) write
\begin{equation}
g(\zeta) =\begin{bmatrix} g_+(\zeta) \\ g_-(\zeta)\end{bmatrix},
\end{equation}
where $g_+(\zeta)\in \mathbb C^{n_1}$ and  $g_-(\zeta)\in \mathbb C^{n_2}$.
Then 
\begin{align*}
  0 = & \tilde{W}_B\begin{bmatrix} ({\mathcal H}S^{-1} g)(1) \\  ({\mathcal H}S^{-1}g)(0) \end{bmatrix}  
     =  \begin{bmatrix}  V_1 & V_2 \end{bmatrix}  \begin{bmatrix} g_+(1) \\  g_-(1)  \end{bmatrix} + \begin{bmatrix}U_1 & U_2 \end{bmatrix}  \begin{bmatrix}  g_+(0) \\  g_-(0)  \end{bmatrix} \\
   =&
   \begin{bmatrix}V_1 & U_2 \end{bmatrix}  \begin{bmatrix} g_+(1) \\  g_-(0)  \end{bmatrix} + \begin{bmatrix} U_1 & V_2 \end{bmatrix}  \begin{bmatrix} g_+(0) \\  g_-(1)  \end{bmatrix} \\
   =&
    \begin{bmatrix} V_1 & U_2 \end{bmatrix} 
    \begin{bmatrix}\Lambda(1)^{-1} & 0 \\ 0 & \Theta(0)^{-1} \end{bmatrix}      \begin{bmatrix} \Lambda (1) g_+(1) \\  \Theta(0) g_-(0)  \end{bmatrix} \\
&+ 
  \begin{bmatrix} U_1 & V_2 \end{bmatrix}  \begin{bmatrix}\Lambda(0)^{-1} & 0 \\ 0 & \Theta(1)^{-1} \end{bmatrix} \begin{bmatrix} \Lambda(0) g_+(0) \\  \Theta(1) g_-(1)  \end{bmatrix}.
\end{align*}
Thus by Proposition  \ref{lemdiag}   the operator $A_S$ as defined in \eqref{eq:10} and \eqref{eq:11} generates a $C_0$-semigroup if and only if 
the matrix 
\[
   K = \begin{bmatrix} V_1 & U_2 \end{bmatrix}
    \begin{bmatrix} \Lambda(1)^{-1} & 0 \\ 0 & \Theta(0)^{-1} \end{bmatrix}
\]
is invertible. 
Since the matrix $ \left[\begin{smallmatrix} \Lambda(1)^{-1} & 0 \\ 0 & \Theta(0)^{-1} \end{smallmatrix}\right]$ is invertible, $A_S$ generates a $C_0$-semigroup if and only if $\begin{bmatrix} V_1 & U_2 \end{bmatrix} $ is invertible. 
Now, $\begin{bmatrix} V_1 & U_2 \end{bmatrix} $ is invertible if and only if for every $f\in \mathbb C^n$ there exists $x\in \mathbb C^{n_1}$ and $y\in \mathbb C^{n_2}$ such that
\begin{eqnarray}
\nonumber
 f &=& \begin{bmatrix} V_1 & U_2 \end{bmatrix} \begin{bmatrix} x \\  y  \end{bmatrix}
= \begin{bmatrix} V_1 & U_2 \end{bmatrix} \begin{bmatrix} x \\  y  \end{bmatrix} +
    \begin{bmatrix} U_1 & V_2 \end{bmatrix}  \begin{bmatrix} 0 \\  0  \end{bmatrix} \\
\nonumber
  &=&   \begin{bmatrix} V_1 & V_2 \end{bmatrix} \begin{bmatrix} x \\  0  \end{bmatrix} +
    \begin{bmatrix} U_1 & U_2 \end{bmatrix}  \begin{bmatrix} 0 \\  y  \end{bmatrix}
\\
\label{eq:12}
   &=&   W_{1}({\mathcal H} S^{-1})(1)\begin{bmatrix} x \\  0  \end{bmatrix}  +   W_{0}({\mathcal H} S^{-1})(0)\begin{bmatrix} 0 \\  y  \end{bmatrix}.
\end{eqnarray}
Referring, to equation (\ref{eqndiagonal}) the columns of $S^{-1}(\zeta)$ are the eigenvectors of $P_1 {\mathcal H}(\zeta)$. 
The eigenvectors corresponding to the  positive eigenvalues forms the first $n_1$ columns. Thus  $S^{-1}(1)\begin{bmatrix} x \\  0  \end{bmatrix}$ is in $Z^+(1)$. Similarly, $S^{-1}(0)\begin{bmatrix} 0\\ y \end{bmatrix} $ is in $Z^-(0)$. 
Thus $\begin{bmatrix} V_1 & U_2 \end{bmatrix} $ is invertible if and only if 
\begin{equation*}
  W_1{\mathcal H}(1)  Z^+(1) \oplus W_0 {\mathcal H}(0) Z^-(0) = {\mathbb C}^n,
\end{equation*}
which concludes the proof.
\hfill${\Box}$

\section{Examples}

The following three examples are provided as illustration of  Theorem \ref{Thmmain}.
\begin{example}
Consider the one-dimensional transport equation on the interval $(0,1)$:
\begin{align*}
\frac{\partial x}{\partial t }(\zeta,t) &= \frac{\partial {\mathcal H}x}{\partial \zeta }(\zeta,t), \qquad x(\zeta, 0) = x_0(\zeta),\\
\begin{bmatrix} w_1 & w_0\end{bmatrix}& \begin{bmatrix}  ({\mathcal H}x)(1,t) \\({\mathcal H}x)(0,t) \end{bmatrix} =0,
\end{align*}
where  $ {\mathcal H}\in C^1[0,1]$ with $ {\mathcal H}(\zeta)>0$ for
every $\zeta\in[0,1]$.

An easy calculation shows $P_1 {\mathcal H}=  {\mathcal H}$ and thus
$Z^+(1)=\mathbb C$ and $ Z^-(0)=\{0\}$. Thus by Theorem \ref{Thmmain} the corresponding operator 
\begin{align*}
A{\mathcal H}x &= \frac{ \partial}{ \partial \zeta} ( {\mathcal H}x),\\
D(A{\mathcal H}) &= \left\{x\in  W^{1,p}(0,1)\mid \begin{bmatrix} w_1 & w_0\end{bmatrix} \begin{bmatrix}  ({\mathcal H}x)(1) \\({\mathcal H}x)(0) \end{bmatrix} =0\right\},
\end{align*}
generates a $C_0$-semigroup on $L^p(0,1)$ if and only if $w_1\not =0$. Further, by Theorem \ref{theo1}, $A{\mathcal H}$ generates a contraction semigroup (unitary $C_0$-group) on $L^2(0,1)$ equipped with
the scalar product $\langle \cdot,  {\mathcal H} \cdot\rangle$ if and only if $w^2_1 \ge w^2_0$ ($w^2_1 = w^2_0$).
\end{example}

\begin{example}\label{ex2}
An (undamped) vibrating string can be modeled by 
  \begin{align}
  \label{eq:7.1.1}
    \frac{\partial^2 w}{\partial t^2} (\zeta,t) = \frac{1}{\rho(\zeta)} \frac{\partial }{\partial \zeta} \left( T(\zeta) \frac{\partial w}{\partial \zeta}(\zeta,t) \right), \qquad t\ge 0, \zeta\in (0,1),
  \end{align}
 where $\zeta\in [0,1]$ is the spatial variable, $w(\zeta,t)$
  is the vertical position of the string at place $\zeta$ and time $t$, $T (\zeta)>0 $ is the Young's modulus
  of the string, and $\rho (\zeta ) >0$ is the mass density, which may vary along the string. We assume that $T$ and $\rho$ are positive and continuously differentiable functions on $[0,1]$. By choosing the state variables 
$x_1= \rho \frac{\partial w}{\partial t}$ (momentum) and $x_2 = \frac{\partial w}{\partial \zeta}$ (strain), the partial differential equation  (\ref{eq:7.1.1}) can equivalently be written as
\begin{align}
  \nonumber
  \frac{\partial }{\partial t} \begin{bmatrix} x_1(\zeta,t) \\ x_2(\zeta,t) \end{bmatrix} &= \begin{bmatrix} 0 & 1 \\ 1 & 0 \end{bmatrix} \frac{\partial }{\partial \zeta}\left( \begin{bmatrix} \frac{1}{\rho(\zeta)} & 0 \\ 0 & T(\zeta) \end{bmatrix}\begin{bmatrix} x_1(\zeta,t) \\ x_2(\zeta,t) \end{bmatrix} \right) \\
 \label{eq:7.2.3} &= P_1  \frac{\partial }{\partial \zeta}\left( {\cal H}(\zeta)\begin{bmatrix} x_1(\zeta,t) \\ x_2(\zeta,t) \end{bmatrix} \right),
\end{align}
where $P_1=\left[\begin{smallmatrix} 0 & 1\\ 1 & 0 \end{smallmatrix}\right]$ and  $ {\cal H}(\zeta)=\left[\begin{smallmatrix}\frac{1}{\rho(\zeta)} & 0\\ 0 & T(\zeta)  \end{smallmatrix}\right]$. 


 The boundary conditions for (\ref{eq:7.2.3}) are 
\[ \begin{bmatrix} W_1 & W_0 \end{bmatrix} \begin{bmatrix}  ({\mathcal H}x)(1,t) \\({\mathcal H}x)(0,t) \end{bmatrix} =0,\]
where $\begin{bmatrix} W_1 & W_0 \end{bmatrix}$ is a $2\times 4$-matrix with rank 2,
or equivalently, the partial differential equation \eqref{eq:7.1.1} is equipped with the boundary conditions
\[  \begin{bmatrix} W_1 & W_0 \end{bmatrix}  \begin{bmatrix}   \rho \frac{\partial w}{\partial t} (1,t)\\\frac{\partial w}{\partial \zeta} (1,t)\\ \rho \frac{\partial w}{\partial t} (0,t) \\\frac{\partial w}{\partial \zeta}(0,t)\end{bmatrix} =0.\]

Defining $\gamma =\sqrt{ T(\zeta)/\rho ( \zeta) }$,  the matrix function $P_1 {\cal H}$ can be factorized as
\[    P_1 {\cal H} = \begin{bmatrix} \gamma & -\gamma \\ \rho^{-1} & \rho^{-1} \end{bmatrix} \begin{bmatrix} \gamma & 0 \\ 0 & -\gamma \end{bmatrix} 
\begin{bmatrix} (2\gamma)^{-1} & \rho/2 \\  (2\gamma)^{-1} & \rho/2 \end{bmatrix}, \]
This implies $Z^+(1) = \mbox{span}\, \left[\begin{smallmatrix} T(1)\\ \gamma(1) \end{smallmatrix}\right]$ and  $Z^-(0) = \mbox{span}\, \left[\begin{smallmatrix} -T(0)\\ \gamma(0) \end{smallmatrix}\right]$. Thus, by Theorem \ref{Thmmain} the corresponding operator 
\begin{align*}
(A{\mathcal H}x)(\zeta) &= \begin{bmatrix} 0 & 1\\ 1 & 0 \end{bmatrix}  \frac{\partial }{\partial \zeta}\left( \begin{bmatrix}\frac{1}{\rho(\zeta)} & 0\\ 0 & T(\zeta)  \end{bmatrix} x(\zeta) \right);\\
D(A{\mathcal H}) &=\left\{x\in  W^{1,p}(0,1;\mathbb C^2)\mid \begin{bmatrix} W_1 & W_0\end{bmatrix} \begin{bmatrix}  ({\mathcal H}x)(1) \\({\mathcal H}x)(0) \end{bmatrix} =0\right\},
\end{align*}
generates a $C_0$-semigroup on $L^p(0,1;\mathbb C^2)$ if and only if 
\[W_1 \begin{bmatrix} \gamma(1) \\ T(1) \end{bmatrix} \oplus W_0 \begin{bmatrix} -\gamma(0) \\ T(0) \end{bmatrix} =\mathbb C^2,\]
or equivalently if the vectors $W_1 \left[\begin{smallmatrix} \gamma(1) \\ T(1) \end{smallmatrix}\right]$ and $W_0 \left[\begin{smallmatrix} -\gamma(0) \\ T(0) \end{smallmatrix}\right]$ are linearly independent.

If $W_1:= I$ and $W_0:= \left[\begin{smallmatrix} -1 & 0\\ 0&1 \end{smallmatrix}\right]$,
then $A{\mathcal H}$ generates a $C_0$-semigroup if and only if the vectors  $ \left[\begin{smallmatrix} \gamma(1) \\ T(1) \end{smallmatrix}\right]$ and $ \left[\begin{smallmatrix} \gamma(0) \\ T(0) \end{smallmatrix}\right]$ are linearly independent. Thus, not only the nature of the boundary conditions but also  Young's modulus
 and the mass density   on the interval $[0,1]$ affect whether or not $A{\mathcal H}$ generates a $C_0$-semigroup.
\end{example}

\begin{example}\label{ex3}
Consider the following network of three  transport equations on the interval $(0,1)$:
\begin{align*}
\frac{\partial x_j}{\partial t }(\zeta,t) &=  \frac{\partial x_j}{\partial \zeta }(\zeta,t), \quad t\ge 0, \, \zeta \in (0,1), \, j=1,2,3,\\
\ x_j(\zeta, 0) &= x_{j,0}(\zeta),\quad \zeta \in (0,1),\, j=1,2,3\\
\begin{bmatrix} 1 & 0 & 0 & 0 & 0 & 0 \\ 0 & 1 & 0 & -1 & 0 & -1 \\ 0 & 0 & 1 & 0 & -1 & 0\end{bmatrix} &\begin{bmatrix}  x_1(1,t) \\ x_2(1,t) \\ x_3(1,t) \\ x_1(0,t) \\ x_2(0,t) \\x_3(0,t) \end{bmatrix} =0, \quad t\ge 0.
\end{align*}
Writing $x= \left[\begin{smallmatrix} x_1 \\ x_2 \\ x_3 \end{smallmatrix}\right]$,  the corresponding operator $A:D(A)\subset L^p(0,1; \mathbb C^3) \rightarrow L^p(0,1; \mathbb C^3)$ is 
\begin{align*}
(Ax)(\zeta) &=\frac{\partial x}{\partial \zeta }(\zeta), \\
D(A) &= \left\{x\in  W^{1,p}(0,1;\mathbb C^3)\mid \begin{bmatrix} 1 & 0 & 0 & 0 & 0 & 0 \\ 0 & 1 & 0 & -1 & 0 & -1 \\ 0 & 0 & 1 & 0 & -1 & 0\end{bmatrix}\begin{bmatrix}  x(1) \\x(0) \end{bmatrix} =0\right\}.
\end{align*}
 In this example  ${\cal H}=P_1=I$ and $P_0=0$  and therefore the assumptions on $S$, $\Lambda$ and $\Theta$ are satisfied.
An easy calculation yields
\begin{align*}
x^*(1)x(1)- x^*(0) x(0) = 2 x_1(0) x_3(0) 
\end{align*}
for every $x\in D(A)$. Theorem \ref{theo1} implies that $A$ does not generate a contraction semigroup on $L^2(0,1;\mathbb C^3)$.

However,  by Theorem \ref{Thmmain}   $A$  generates a $C_0$-semigroup on $L^p(0,1;\mathbb C^3)$ for $1\le p< \infty$:  In this example,  $Z^+(\zeta)= \mathbb C^3$ , $Z^-(\zeta)= \{0\}$, $W_1= I$ and $W_0 =\left[\begin{smallmatrix}  0 & 0 & 0 \\  -1 & 0 & -1 \\  0 & -1 & 0\end{smallmatrix}\right]$. Thus,
\begin{align*}
W_1 Z^+(1) \oplus W_0 Z^-(0) = \mathbb C^3.
\end{align*}
Finally, \cite[Corollary 2.1.6]{S} implies that  $A$ generates a contraction semigroup on $L^1(0,1;\mathbb C^3)$.
  
   Summarizing, $A$ generates a generates a $C_0$-semigroup on $L^p(0,1;\mathbb C^3)$ for $1\le p< \infty$ and in fact a  contraction semigroup on  $L^1(0,1;\mathbb C^3)$ but it does not generate a contraction  semigroup on  $L^2(0,1;\mathbb C^3)$.
\end{example}

\section*{Acknowledgement}

The authors gratefully acknowledge support from the DFG (Grant
JA 735/9-1) and NWO (Grant DN 63-261)


\begin{thebibliography}{9}

\bibitem{AJ}
\newblock B. Augner and B. Jacob,
 \newblock Stability and stabilization of infinite-dimen\-sional
linear port-Hamiltonian systems, {\em Evolution Equations and Control Theory}, \textbf{3}(2) (2014), 207--229.

\bibitem{E}
\newblock K.-J.~Engel,
\newblock Generator property and stability for generalized difference operators,
{\em Journal of Evolution Equations}, \textbf{13}(2) (2013), 311--334.


\bibitem{JZ}
 \newblock B.~Jacob and H.J. Zwart,
 \newblock {\em Linear Port-Hamiltonian Systems on Infinite-dimensional Spaces},
 \newblock Operator Theory: Advances and Applications, \textbf{223} (2012),
Birkh\"auser, Basel.

\bibitem{GZM}
 \newblock Y. Le~Gorrec, H. Zwart and B. Maschke,
 \newblock Dirac structures and boundary control systems associated with skew-symmetric differential operators,
 \newblock {\em SIAM J. Control Optim.}, \textbf{44} (2005), 1864--1892.

\bibitem{S}
\newblock E.~Sikolya,
\newblock Semigroups for flows in networks,
\newblock Ph.D thesis, University of T\"ubingen, 2004.
 

\bibitem{VanDerSchaftMaschke_2002} 
 \newblock A.J. van der Schaft and B.M. Maschke,
 \newblock Hamiltonian formulation of distributed parameter systems with boundary energy flow,
 \newblock {\em J. Geom. Phys.}, \textbf{42} (2002), 166--174.


\bibitem{Villegas_2007}
 \newblock J.A. Villegas,
 \newblock {\em A port-Hamiltonian Approach to Distributed Parameter Systems},
 \newblock Ph.D thesis, Universiteit Twente in Enschede, 2007. Available from: \url{http://doc.utwente.nl/57842/1/thesis_Villegas.pdf}.

\bibitem{ZGMV}
 \newblock H. Zwart, Y. Le Gorrec, B. Maschke and J.
              Villegas,
 \newblock Well-posedness and regularity of hyperbolic boundary control
              systems on a one-dimensional spatial domain,
 \newblock {\em ESAIM Control Optim. Calc. Var.}, \textbf{16}(4) (2010), 1077-1093.
 
\bibitem{K}
\newblock T. Kato, 
\newblock {\em Perturbation Theory for Linear Operators},
\newblock Springer-Verlag, Berlin, 1995.
    
   
\end{thebibliography}
\end{document}